\documentstyle[twoside, 12pt, amstex]{amsart}
\makeatletter

\begin{document}

\newcommand{\e}{\epsilon}
\newcommand{\w}{{\bold w}}
\newcommand{\y}{{\bold y}}
\newcommand{\z}{{\bold z}}
\newcommand{\x}{{\bold x}}
\newcommand{\N}{{\mathbb N}}
\newcommand{\Z}{{\mathbb Z}}
\newcommand{\F}{{\bold F}}
\newcommand{\R}{{\mathbb R}}
\newcommand{\BR}{{\mathbb R}}
\newcommand{\Q}{{\mathbb Q}}
\newcommand{\Ql}{{\mathbb Q_{\ell}}}
\newcommand{\C}{{\mathbb C}}
\newcommand{\K}{{\mathbb K}}
\newcommand{\G}{{\mathbf G}}
\newcommand{\BA}{{\mathbb A}}
\newcommand{\BC}{{\mathbb C}}
\newcommand{\BP}{{\mathbb P}}
\newcommand{\BQ}{{\mathbb Q}}
\newcommand{\BZ}{{\mathbb Z}}
\newcommand{\cF}{{\mathcal F}}
\newcommand{\cL}{{\mathcal L}}
\newcommand{\cR}{{\mathcal R}}
\newcommand{\cO}{{\mathcal O}}
\newcommand{\cX}{{\mathcal X}}
\newcommand{\cH}{{\mathcal H}}
\newcommand{\cM}{{\mathcal M}}
\newcommand{\sB}{{\sf B}}
\newcommand{\cC}{{mathacal C}}
\newcommand{\cT}{{\mathcal T}}
\newcommand{\cI}{{\mathcal I}}
\newcommand{\cS}{{\mathcal S}}
\newcommand{\sE}{{\sf E}}

\newcommand{\gp}{{\mathfrak p}}

\newcommand{\sA}{{\sf A}}
\newcommand{\ga}{{\sf a}}
\newcommand{\es}{{\sf s}}
\newcommand{\m}{{\bold m}}
\newcommand{\bS}{{\bold S}}
\newcommand{\Fp}{{\mathbb F}_p}
\newcommand{\Fq}{{\mathbb F}_q}
\newcommand{\tr}{\operatorname{tr}}
\newcommand{\Id}{\operatorname{Id}}

\newcommand{\ihra}{\stackrel{i}
{\hookrightarrow}}
\newcommand\rank{\mathop{\rm rank}\nolimits}
\newcommand\im{\mathop{\rm Im}\nolimits}
\newcommand\Li{\mathop{\rm Li}\nolimits}
\newcommand\NS{\mathop{\rm NS}\nolimits}
\newcommand\Hom{\mathop{\rm Hom}\nolimits}
\newcommand\Pic{\mathop{\rm Pic}\nolimits}
\newcommand\Spec{\mathop{\rm Spec}\nolimits}
\newcommand\Hilb{\mathop{\rm Hilb}\nolimits}
\newcommand{\length}{\mathop{\rm length}
\nolimits}

\newcommand\lra{\longrightarrow}
\newcommand\ra{\rightarrow}
\newcommand\cJ{{\mathcal J}}
\newcommand\JG{J_{\Gamma}}
\newcommand{\wvskp}{\vspace{1cm}}
\newcommand{\vskp}{\vspace{5mm}}
\newcommand{\nvskp}{\vspace{1mm}}
\newcommand{\nid}{\noindent}
\newcommand{\new}{\nvskp \nid}
\newtheorem{Assumption}{Assumption}[section]
\newtheorem{Theorem}{Theorem}[section]
\newtheorem{Lemma}{Lemma}[section]
\newtheorem{Remark}{Remark}[section]
\newtheorem{Corollary}{Corollary}[section]
\newtheorem{Conjecture}{Conjecture}[section]
\newtheorem{Proposition}{Proposition}[section]
\newtheorem{Example}{Example}[section]
\newtheorem{Definition}{Definition}[section]
\newtheorem{Question}{Question}[section]
\newtheorem{Problem}{Problem}[section]

\renewcommand{\thesubsection}{\it}

\title{The modularity of certain non-rigid
Calabi--Yau threefolds}

\author{Ron Livn\'e}
\thanks{Ron Livn\'e was partially supported
by a Israel-USA BSF Research Grant.}
\address{Institute of Mathematics, Hebrew
University of Jerusalem, Givat-Ram, Jerusalem
91904, Israel}
\email{rlivne@math.huji.ac.il}
\author{Noriko Yui}
\thanks{Noriko Yui was partially supported by a
Discovery Grant from NSERC, Canada.}
\address{Department of Mathematics, Queen's
University, Kingston. Ontario Canada K7L 3N6}
\email{yui@mast.queensu.ca;
yui@fields.utoronto.ca}


\date{Version of June 16, 2005}

\begin{abstract}
Let $X$ be a  Calabi--Yau threefold fibred over
$\BP^1$ by non-constant semi-stable K3 surfaces
and reaching the Arakelov--Yau bound. In [STZ],
X. Sun, Sh.-L. Tan, and K. Zuo  proved that $X$
is modular in a certain sense. In particular, the
base curve is a modular curve. In their result
they distinguish the rigid and the non-rigid
cases. In [SY] and [Y] rigid examples were
constructed. In this paper we construct explicit
examples in non-rigid cases. Moreover, we prove
for our threefolds that the ``interesting'' part
of their $L$\/-series is attached to an
automorphic form, and hence that they are modular
in yet another sense.
\end{abstract}

\maketitle
\section{Introduction}
Let $X$ be an algebraic threefold and let $f:X\to
\BP^1$ be a non-isotrivial morphism whose fibers
are semi-stable K3 surfaces. Let $S\subset \BP^1$
be the finite set of points above which $f$ is
non-smooth, and assume that the monodromy at each
point of $S$ is non-trivial.  Jost and Zuo [JZ]
proved the Arakelov--Yau type inequality:
\[\mbox{deg} f_* \omega_{X/\BP^1}\leq \mbox{deg}\,
\Omega_{\BP^1}^1 (\mbox{log} S).\]

Let $\Delta\subset X$ be the pull-back of $S$.
Let $\omega_{X/{\BP}^1}$ be the canonical sheaf.
The Kodaira--Spencer maps $\theta^{2,0}$ and
$\theta^{1,1}$ are maps of sheaves fitting into
the following diagram:
 \[f_*\Omega^2_{X/{\BP}^1}(\mbox{log}\,\Delta)
\overset{\theta^{2,0}}\to
R^1f_*\Omega^1_{X/{\BP}^1}(\mbox{log}\,\Delta)
\otimes \Omega^1_{{\BP}^1}(\mbox{log}\,S)
\overset{\theta^{1,1}}\to R^2\,
f_*\cO_{X/{\BP}^1}\otimes\Omega^1_{{\BP}^1}
(\mbox{log}\,S)^{{\otimes}2}.\]
 The iterated Kodaira--Spencer map of
$f$ is defined to be the map
$\theta^{1,1}\theta^{2,0}$.
\smallskip

It is known (see [STZ]) that when the (iterated)
Kodaira--Spencer map is $0$, one actually has the
stronger inequality
\[\mbox{deg} f_* \omega_{X/\BP^1}\leq
\frac{1}{2}\mbox{deg}\,
\Omega_{\BP^1}^1 (\mbox{log} S).\]

Assume from now on that $X$ is a Calabi--Yau
threefold. Then the triviality of the canonical
bundle implies that $\mbox{deg}\, f_*
\omega_{X/\BP^1}=2$ (see (\ref{deg}) below).

Recently X. Sun, S-L. Tan and K. Zuo [STZ]
considered Calabi--Yau threefolds for which the
Arakelov--Yau inequality becomes equality. Thus
$S$ consists of $6$ points when the
Kodaira--Spencer map is $0$ and of $4$ points
otherwise.

As a consequence of the main theorem of [STZ],
the following results were obtained.
\medskip

{\bf Theorem 1} (Corollary 0.4 in
[STZ]): {\sl (i) If the iterated Kodaira--Spencer
map $\theta^{1,1}\theta^{2,0}$ of $f$ is
non-zero, then $X$ is rigid (i.e., $h^{2,1}=0$)
and birational to the Nikulin--Kummer
construction of a symmetric square of a family of
elliptic curves $f: E\to \BP^1$.  After passing
to a double cover $E^{\prime}\to E$ (if
necessary), the family $g^{\prime}: E^{\prime}\to
\BP^1$ is one of the six modular families of
elliptic curves on the Beauville's list ([B]).
\smallskip

(ii) If the iterated Kodaira--Spencer map
$\theta^{1,1}\theta^{2,0}$ of $f$ is zero, then
$X$ is a non-rigid Calabi--Yau threefold (i.e.,
$h^{2,1}\neq 0$), the general fibers have Picard
number at least $18$, and $\BP^1\setminus
S\simeq \mathfrak H/\Gamma$ where $\Gamma$ is a
congruence subgroup of $PSL(2,\Z)$ of index $24$.}
\medskip

{\bf Remark}: The base curve
$\BP^1\setminus S$ is a modular variety of genus
zero, i.e.,  $\mathfrak H/\Gamma$ where $\Gamma$
is a torsion-free genus zero congruence subgroup
of $PSL(2,\BZ)$ of index $12$ in case (i), and of
index $24$ in case (ii). In the paper of Sun, Tan
and Zuo [STZ], the word ``modularity'' refers to
this fact.
\medskip

The third cohomology of each of the six rigid
Calabi--Yau threefolds in Theorem 1 (i) arises
from a weight $4$ modular form. In the articles
of Saito and Yui [SY] and of Verrill in Yui [Y],
these forms were explicitly determined. Saito and Yui
use geometric structures; while Verrill uses point
counting method, to obtain the results. 

More precisely, the following was proved for the
natural models over $\BQ$ of these six rigid
threefolds:
\medskip

{\bf Theorem 2} (Saito and Yui [SY] and
Verrill in Yui [Y]): {\sl For each of the six
rigid Calabi--Yau threefold over $\Q$, the
L-series of the third cohomology coincides with
the L-series arising from the cusp form of weight
$4$ of one variable on the modular group in the
Beauville's list. Beauville's list and the
corresponding cusp forms are given in Table 1.}
\smallskip
\[
\begin{array}{|c|c|c|} \hline \hline
\mbox{Group} & \mbox{Number of components} &
\mbox{Cusp forms}    \\
\Gamma       & \mbox{of singular fibers}   &
\mbox{of weight $4$}  \\ \hline
\Gamma(3)   & 3,3,3,3                &
\eta(q^3)^8   \\ \hline
\Gamma_1(4)\cap \Gamma(2) & 4,4,2,2  &
\eta(q^2)^4\eta(q^4)^4 \\ \hline
\Gamma_1(5) & 5,5,1,1                &
\eta(q)^4\eta(q^5)^4   \\ \hline
\Gamma_1(6) & 6,3,2,1              &
\eta(q)^2\eta(q^2)^2\eta(q^3)^2\eta(q^6)^2 \\
\hline
\Gamma_0(8)\cap\Gamma_1(4) & 8,2,1,1&
\eta(q^4)^{16}\eta(q^2)^{-4}\eta(q^8)^{-4} \\
\hline
\Gamma_0(9)\cap\Gamma_1(3) & 9,1,1,1 &
\eta(q^3)^8      \\ \hline \hline
\end{array}
\]
\centerline{Table 1: Rigid Calabi--Yau
threefolds and cusp forms}
\smallskip

\noindent Here $\eta(q)$ denotes the Dedekind eta-function:
$\eta(q)=q^{1/24}\prod_{n\geq 1}(1-q^n)$ with $q=e^{2\pi i\tau}$.

\pagebreak

It might be helpful to recall the six rigid Calabi--Yau
threefolds in Theorem 2. These six rigid Calabi--Yau threefolds 
are obtained (by Schoen [S]; 
see also Beauville [B]) as the self-fiber products of stable families 
of elliptic curves admitting exactly
four singular fibers. The base curve is a rational modular curve
and correspond to the torsion-free genus zero
congruence subgroups 
$\Gamma$ of $PSL(2,\BZ)$ in Table 1.   Note that the
$4$-tuples of natural numbers appearing in the
second column add up to $12$, which is the index of
the modular group $\Gamma$ in $PSL(2,\BZ)$.
\medskip

In [STZ] the authors indicate one example for the
non-rigid extremal case. It is related to
$\Gamma(4)$, which is a torsion-free congruence
subgroup of genus $0$ and index $24$ in
$PSL(2,\BZ)$. The list of torsion-free congruence
subgroups of genus $0$ and index 24 in
$PSL(2,\BZ)$ is known (see Sebbar [Se], and Table
2 below). In this paper we will show that most of
them  give rise to a similar collection of
examples. In each of these cases we will compute
the interesting part of the $L$\/-series of the
third cohomology of an appropriate natural  model
over $\BQ$ in terms of automorphic forms.

This paper is organized as follows. In Section 2, 
we use work of Sebbar [Se] to determine the
groups $\Gamma$ corresponding to case (ii) of
Theorem 2. These are subgroups of $PSL(2,\BZ)$,
and associated to each $PSL(2,\BZ)$\/-conjugacy
class there is a natural elliptic fibration over
the base curve, defined over $\BQ$. The total
spaces of these fibrations are elliptic modular
surfaces in the sense of Shioda [Sh1]. Moreover
each is an extremal K3 surface (namely their
Picard number is 20, the maximum possible). We
explain the relation between the motive of their
transcendental cycles, and specific CM forms of
weight $3$ using a result of Livn\'e on
orthogonal rank 2 motives in [L2].

Section 3 contains our main results: we
construct our examples, verify the required
properties, and analyze the interesting part of
their cohomology. (See the final Remarks 8 (2) for the
other parts.) Then in Section 4 we give explicit
formulas for the weight $3$ cusp forms and
defining equations for the elliptic fibrations of
Section 2.

The paper is supplemented by the article of Hulek and Verrill
[HV] where they treat Kummar varieties, one of which is  
the case associated to the modular group $\Gamma_1(7)$. 
This case differs from the examples considered in our
paper with the main difference being the fact that the $2$-torsion
points do not decompose into four sections, leading to non-semi-stable
fibrations.  But it still gives rise to a Calabi--Yau threefold (Theorem 2.2
of Hulek and Verrill [HV]), and the modularity question 
can still be considered, and this is exactly what Hulek
and Verrill deals with in the supplement [HV] to this article.

\section{Extremal congruence K3 surfaces}

The torsion-free genus zero congruence subgroups
of $PSL(2,\Z)$ of index $24$ were classified in
Sebbar [Se]. There are precisely nine conjugacy
classes of such congruence subgroups.
\medskip

The second column in the following Table 2 gives
the complete list of the torsion-free congruence
subgroups of $PSL(2,\Z)$ of index $24$ up to
conjugacy. Each has precisely 6 cusps. The third
column in the table gives the widths of these
cusps.

\[
\begin{array}{|c|c|c|} \hline\hline
\#& \mbox{The group $\Gamma$} &
\mbox{Widths of the cusps}  \\
\hline
1& \Gamma(4) & 4,4,4,4,4,4  \\ \hline
2& \Gamma_0(3)\cap \Gamma(2) & 6,6,6,2,2,2  \\
\hline
3& \Gamma_1(7) & 7,7,7,1,1,1  \\ \hline
4& \Gamma_1(8) & 8,8,4,2,1,1  \\ \hline
5& \Gamma_0(8)\cap \Gamma(2) & 8,8,2,2,2,2  \\
\hline
6& \Gamma_1(8;4,1,2) & 8,4,4,4,2,2  \\ \hline
7& \Gamma_0(12) & 12,4,3,3,1,1  \\ \hline
8& \Gamma_0(16) & 16,4,1,1,1,1  \\ \hline
9& \Gamma_1(16;16,2,2) & 16,2,2,2,1,1  \\
\hline \hline
\end{array}
\]
\centerline{Table 2: Torsion-free congruence
subgroups of index 24.}
\smallskip

\noindent Here
 \[\Gamma_1(8;4,1,2):=\{\pm
\begin{pmatrix} 1+4a & 2b \\ 4c & 1+4d
\end{pmatrix},\, a\equiv c\pmod 2\}\] and
\[\Gamma_1(16;16,2,2):=\{\pm \begin{pmatrix}
1+4a & b \\
8c & 1+4d \end{pmatrix}, \, a\equiv c\pmod 2\}.\]
\medskip 
{\bf Remark 3}: {\sl If we are
interested in conjugacy as Fuchsian groups (in
$PSL(2,\BR)$), Examples \#1,\#5, and \#8 are conjugate
(use the matrix $\begin{pmatrix} 1 & 0
\\ 0 & 2
\end{pmatrix})$, Examples \#2 and \#7 are conjugate
(use the same matrix), and Examples \#4, \#6, and \#9
are conjugate (use the same matrix as well as
$\begin{pmatrix} 1 & 1 \\ 0 & 1 \end{pmatrix})$.}
\medskip

{\bf Proposition 4}: {\sl Let $\Gamma$
be one of the congruence subgroups in Table 2.
Then $\Gamma$ has an explicit congruence lift
$\tilde{\Gamma}$ to $SL(2,\BZ)$ with the
following properties:

(1) $\tilde{\Gamma}$ has no elliptic elements.
In particular $-{\rm Id}$ is not in $\tilde{\Gamma}$.

(2) $\tilde{\Gamma}$ contains no element of trace
$-2$.}

{\bf Proof}: We let $\tilde{\Gamma}$ be
the subgroup of $SL(2,\BZ)$ consisting of the
elements above $\Gamma$. 
Indeed, the lifts $\tilde{\Gamma}$ are sometimes the same as the groups 
$\Gamma$ themselves. In fact, for the cases
\#1,\,\#2,\,\#3,\, \#4 and \#5, the lifts are the same and
respectively given by: $\Gamma(4)$,\,
$\Gamma_0(3)\cap\Gamma(2)$,\,$\Gamma_1(7)$,\,$\Gamma_1(8)$ and
$\Gamma_0(8)\cap\Gamma(2)$. 
In the cases \#6,\,\#7,\,\#8 and \#9, the lifts are
not unique, and we choose respectively
the following lifts: $\Gamma_1^{\prime}(8;4,1,2)$,\, $\Gamma_0(12)\cap\Gamma_1(3)$, $\Gamma_0(16)\cap\Gamma_1(4)$ and $\Gamma_1^{\prime}(16;16,2,2)$.
Here $\Gamma_1^{\prime}(8;4,1,2)$ and $\Gamma_1^{\prime}(16;16,2,2)$
are defined in the same way as $\Gamma_1(8;4,1,2)$ and $\Gamma_1(16;16,2,2)$
but without the $\pm$. 
Note that the widths of the cusps are not affected by
taking a lift as $-{\rm Id}$ is the only difference.

(We should remark that the lifts are not unique; for instacne, in the case
\#7, there are four lifts, but only one has no elements of
trace $-2$, which is the one given above.)  \qed
\medskip

In fact, Proposition 4 has also been obtained by A. Sebbar in his 
unpublished note. 

Proposition 4 will pave a way to the definition of 
elliptic modular surfaces, which we will discuss next.  
\medskip

{\bf Elliptic Modular Surfaces}: In
[Sh1] Shioda has shown how to associate to any
subgroup $G$ of $SL(2,\BZ)$ of finite index and
not containing $-$Id an elliptic fibration
$E(G)$, called the elliptic modular surface
associated to $G$, over the modular curve
$X(G)=\overline{G\backslash{\frak{H}}}$.

(Shioda's construction requires that modular groups
ought to be a subgroup of $SL(2,\BZ)$ (rather than a subgroup
of $PSL(2,\BZ)$) that contains
no element of order $2$. This is a reason we consider a
lift $\tilde{\Gamma}$ of $\Gamma$ in our discussion.)
\medskip

{\bf Remark 5\ }: {\sl It follows from
Kodaira's theory that when $G$ contains no
elliptic elements and no elements of trace $-$2,
all the singular fibers are above the cusps and
are of type $I_n$, where $n$ is the width of the
cusp. On the other hand, elements of trace $-2$
give rise to $I_n^*$\/-fibers above the cusps.}
\medskip

For the $\tilde{\Gamma}$\/'s of Proposition 4
the modular curve $X(\tilde{\Gamma})$ has genus
$0$, and, since the sum of the widths of the
cusps in Table 2 is always 24, each
$E(\tilde{\Gamma})$ is an extremal K3 surface.
The space $S_3(\tilde{\Gamma})$ of cusp forms of
weight 3 for $\tilde{\Gamma}$ is therefore
one-dimensional. Up to a square, the
discriminant of the intersection form on the
rank 2 motive of the transcendental cycles $T: =
T(E(\tilde{\Gamma})) =
          H^2(E(\tilde{\Gamma}),\BQ)/{\rm NS}
          (E(\tilde{\Gamma}))$
is given by

\begin{equation}
\label{delta}
\delta = \delta_k = -1,-3,-7,-2,-1,-2,-3,-1,-2
\end{equation}
in cases $k = 1,\dots,9$ respectively. To see
this one computes the discriminant of the
(known) lattice $NS(E(\tilde{\Gamma}))$ and
passes to the orthogonal complement in
$H^2(E(\tilde{\Gamma}),\BQ)$. For details, see
e. g. Besser-Livn\'e [BS]. By [L2] it follows
that the normalized newform $g_{3,\Gamma}$
generating $S_3(\tilde{\Gamma})$ has CM by
$\BQ(\sqrt{\delta}\,)$.

To each of our 9 examples there is a naturally
associated moduli problem of classifying
(generalized) elliptic curves with a given level
structure (Katz and Mazur [KM]). Each of these
moduli problems refines the respective moduli
problem $Y_1(M)$, of classifying elliptic curves
with a point of order $M$, where $M$ is as
above. Since $M\geq 4$, these moduli problems
are all represented by universal families
$E(\tilde{\Gamma})/X(\tilde{\Gamma})/\BZ[1/M]$
(see Katz--Mazur [KM]). The geometric fibers are
geometrically connected in all these examples,
and their compactified fibers over $\BC$ are the
corresponding elliptic modular surfaces above.

We shall now compute the $L$\/-series $L(T,s)$ of
the transcendental cycles. By the
Eichler--Shimura Isomorphism, this is the
parabolic cohomology
\[ \tilde{H}: =
\tilde{H}{\vphantom{H}}^1(
X(\tilde{\Gamma})\times_{\BZ[1/M]}\overline{\BQ},
R^1(E(\tilde{\Gamma})\rightarrow
X(\tilde{\Gamma}))).\] Moreover, Deligne proved
([D]) the Eichler--Shimura congruence relation
\[\text{Frob}_p + \text{Frob}'_p = T_p
\qquad\text{for any $p\not|M$},\]
 where $T_p$ is the $p$\/-th Hecke operator on
$S_3(\tilde{\Gamma})$. This is the same as the
$p$\/-th Fourier coefficient of the normalized
newform $g_{3,\Gamma}$. Summarizing, we proved
\[L(T(X(\tilde{\Gamma}),s) = L(g_{3,\Gamma},s).\]

Explicit Weierstrass equations for the elliptic
fibrations
$E(\tilde{\Gamma})/X(\tilde{\Gamma})/\BZ[1/M]$
will be given in Section 4 below.
\medskip

{\bf Remarks}: {\sl The list in Table 2
exhausts all of the families of semi-stable
elliptic $K3$ surfaces with exactly six singular
fibers, which correspond to torsion-free genus
zero {\it congruence} subgroups of $PSL_2(\Z)$ of
index $24$. The $6$-tuples of natural numbers
appearing in the third column add up to $24$.
Therefore, the number of such $6$-tuples is a
priori finite. That this list is complete was
proved by Sebbar [Se].
\smallskip

(2) Miranda and Persson [MP] studied all possible
configurations of $I_n$ fibers on elliptic $K3$
surfaces.  In the case of exactly six singular
fibers, they obtained $112$ possible
configurations including the above nine cases.
All these $K3$ surfaces have the maximal possible
Picard number $20$. It should be emphasized that
the exactly nine configurations correspond to
genus zero congruence subgroups of $PSL_2(\Z)$ of
index $24$.
\smallskip

(3) The theory of Miranda and Persson had been
extended to prove the uniqueness (over $\BC$)  of
K3 surfaces having each of these types of
singular fibers by Artal-Bartolo, Tokunaga and
Zhang [BHZ]. Confer also the article of Shimada and Zhang [SZ] for
a useful table of extremal elliptic $K3$ surfaces.}
\medskip

\section{The non-rigid examples}
Let $Y=E(\tilde{\Gamma})$ be one of the K3
surfaces of the previous Section, and let $g_Y =
g_{3,\Gamma}$ denote the corresponding cusp form
of level $3$. If $\tilde{\Gamma}$ contains
$\Gamma_1(M)$ (in Table 2 this happens in cases
\#3, \#4, \#7, \#8, and \#9), and $M = M_Y$ is
the maximal possible, then $M$ is the {\em level}
of $g_Y$, and the   {\em newtype} of $g_Y$ is the
Dirichlet character
\begin{equation}
\label{new}
\epsilon = \epsilon_Y,
\end{equation}
of conductor $M$, so that $g_Y$ is in
$S_3(\Gamma_0(M), \epsilon_Y)$. Notice that
$\epsilon$ is odd (namely $\epsilon(-1) = -1$).
Moreover, since the coefficients of $g_Y$ are
integers, $\epsilon$ must be quadratic. (We will
determine $\epsilon$ for our examples in
Proposition~10 below.)

Let $E$ be an elliptic curve.  We view the
product $Y\times E$ as a family of abelian
surfaces over $X(\tilde{\Gamma})$. The fiber $A_t
= A_{\Gamma,t}$ over each point $t\in
X(\tilde{\Gamma})$ is the product of the fiber
$E_{\Gamma,t}$ of $E(\tilde{\Gamma})$ with $E$.
Then we have the following
\medskip

{\bf Theorem 6}: {\sl (1) The product
$Y\times E$ has the Hodge  numbers
\[h^{0,3}(Y\times E)=1,\, h^{1,0}(Y\times E) =
1\,\,\mbox{and}\,\, B_3(Y\times E)=44\]
 (so that
$Y\times E$ is not a Calabi--Yau threefold).

(2) The motive $T(Y\times E)=T(Y)\times H^1(E)$
is a submotive of $H^3(Y\times E)$. If $E$ and
$Y$ are defined over $\BQ$, this submotive is
modular, in the sense that its $L$\/-series is
associated to a cusp form $g_Y$ on
$GL(2,\BQ(\sqrt{\delta}\,))$.
\newline Let $g_E$ be the cusp form of
weight $2$ associated to $E$ by Wiles et. al.
([W]). Let $A(p)$ (respectively $B(p)$\/) be the
$p$\/th Fourier coefficient of $g_E$
(respectively of $g_Y$\/), and let $\epsilon_Y$
be the newtype character of $g_Y$ (see Section
3). Then for any good prime $p$, the local
Euler factor $L_p(s)$ of the $L$\/-series
$L(T(Y\times E),s) = L(g_E\otimes g_Y,s)$ is
\[1 - A(p)B(p)p^{-s} + (B(p)^2+\epsilon_Y(p)pA(p)^2-2p^2
 \epsilon_Y(p))p^{1-2s}
  - A(p)B(p)\epsilon_Y(p)p^{3-3s} + p^{6-4s}.\]
}
\medskip

{\bf Proof}: The
statements about the Hodge and Betti numbers
follow from the K\"unneth formula. Since $T(Y)$
is a factor of $H^2(Y)$, it follows again from
the K\"unneth formula that $T(Y)\times H^1(E)$
is a factor of $H^3(Y\times E)$.
\medskip

For the second part, we know that $g_Y$ is a CM
form. Hence it is induced from an algebraic Hecke
character $\chi = \chi_Y$ of the imaginary
quadratic field $F = K_i$. Let $\chi_G$ be the
compatible system of $1$\/-dimensional
$\ell$\/-adic representations of $G_F =
\text{Gal}(\overline\Q/F)$ corresponding to
$\chi$. Then the $2$\/-dimensional Galois
representation associated to $T(Y)$ is
ind$_{G_F}^{G_\Q}\chi_G$. Hence we obtain the
$4$\/-dimensional Galois representation
\[\rho_E \otimes \text{ind}_{G_F}^{G_\Q} \chi_G
\simeq \text{ind}_{G_F}^{G_\Q}(\chi_G \otimes
\text{Res}^{G_\Q}_{G_F}\rho_E),\]
 where $\rho_E$ is the Galois representation
associated to $H^1(E)$. The operation of
restricting $\rho_E$ to $G_F$ and of twisting by
characters have automorphic analogs. Let $\pi_E$
be the automorphic representation associated to
$E$. Then $\pi' = \chi\otimes$ Res$^\Q_F\pi_E$
makes sense as an automorphic cuspidal
irreducible representation of $GL(2,F)$, and we
have the characterizing relationship
\[L(\pi',s) = L(\text{ind}_F^\Q\pi',s) =
L(\pi_E \otimes \text{ind}_F^\Q \chi,s) =
L(g_E\otimes g_Y,s).\]

For the last part, write the $p$\/th Euler
factors of $g_E$ and $g_Y$ respectively as
\[(1-\alpha_p p^{-s})(1-\alpha'_p p^{-s}) =
1-A(p) p^{-s} + p^{1-2s} \qquad\text{and}\]
\[(1-\beta_p p^{-s})(1-\beta'_p p^{-s}) =
1 - B(p) p^{-s} +\epsilon_Y(p) p^{2-2s}.\]
 Then the Euler factor $L_p$ is defined as
\[L_p(s) =(1-\alpha_p\beta_p p^{-s})
(1-\alpha'_p\beta_p p^{-s})
          (1-\alpha_p\beta_p' p^{-s})
          (1-\alpha_p'\beta_p' p^{-s}),\]
and the claim follows by a direct calculation.
\qed
\medskip

{\bf Remark}: {\sl If a K3 surface has the
Picard number $19$ or $18$, the modularity
question for the product $Y\times E$ is still
open. However, if the Picard number is 19, one
knows at least that the rank $3$ motive $T(Y)$
of the transcendental cycles is self dual
orthogonal via the cup product. (For explicit
examples of K3 surfaces with Picard number 19,
see e. g. Besser and Livn\'e [BL].) Thus one can
use a result of Tate to lift each $\ell$\/-adic
representation to the associated spin cover,
which is the multiplicative group of some
quaternion algebra over $\Q$. If the spin
representation is modular (which should always
be the case), then it is associated to a cusp
form $h$ of weight $2$ on $GL(2)$, so that
Symm$^2 h$ realizes $T(Y)$. Let $g_E$ be again
the weight $2$ cusp form associated with $E$. It
follows, by work of Gelbart and Jacquet, that
$T(Y)$ is realized by an automorphic
representation on $GL(3,\Q)$. Hence, by the work
of Kim and Shahidi ([KS]), $T(Y)\times E$ is
realized by an automorphic form on $GL(6,\BQ)$.
In particular, $L(\text{Symm}^2h \otimes g_E,s)$
has the expected analytic properties.}
\medskip

To construct our promised examples, let $X=
X_\Gamma\ra X(\tilde{\Gamma})$ be the associated
Kummer family, in which we divide each fiber
$A_t$ of $Y(\tilde{\Gamma})\times E$ by $\pm 1$
and then blow up the locus of points of order
$2$. We now have the following
\medskip

{\bf Theorem 7\ }: {\sl In the Examples
\#1, \#2, \#5, and \#6 of Table 2 the resulting $X$ is a
smooth Calabi--Yau threefold.  It is non-rigid,
and the given fibration $f:X\ra
X(\tilde{\Gamma})$ is semi-stable, with
vanishing (iterated) Kodaira--Spencer mapping.
We have
\[\deg{f_*\omega_{X/\BP^1}} = 2 =
\frac{1}{2}\,\mbox{deg}\,\Omega^1_{\BP^1}
(\mbox{log}\,S),\]
 In other words, $X$ reaches the (stronger)
 Arakelov--Yau  bound. }
\medskip

{\bf Remark}: {\sl For the first  case
in Table 2 ($\tilde{\Gamma} = \Gamma(4)$) this
example is indicated in [STZ].}
\medskip

{\bf Proof}:  The
Examples we chose are those in which
$\tilde{\Gamma}$ is a subgroup of $\Gamma(2)$.
(This is because otherwise, the points of order $2$
of $X(\Gamma)$ coincide (over the cusps).) Thus the $2$\/-torsion points (of $E_t$ and
hence of $A_t$) are distinct for {\em all} $t\in
X(\tilde{\Gamma})$. It follows that the locus
$A[2]$ of $2$\/-torsion points is smooth, and
hence so is the blow-up $X$. We have $H^i(X) =
H^i(Y(\tilde{\Gamma})\times E)^{<\pm 1>}$. But
$\pm 1$ acts as $\pm 1$ on both the non-trivial
holomorphic $1$\/-form $\omega_1$ of $E$ and on the
non-trivial holomorphic $2$\/-form $\omega_2$ of
$Y(\tilde{\Gamma})$. Hence $\omega_1\wedge
\omega_2$ descends to a holomorphic $3$\/-form
$\omega_3$ on $X$. Its divisor can only be
supported on the proper transform
${\mathcal{F}}$ of $A[2]$; however $\mathcal{F}$
intersects each fiber $f^{-1}(t)$ in sixteen
$(-2)$\/-curves, which do not contribute to the
canonical class, so that $\omega_3$ is indeed
nowhere-vanishing. The K\"unneth formula gives
that
\[H^1(X) = H^1(Y)^{<\pm 1>} =
H^1(E)^{<\pm 1> } = 0, \qquad\text{and}\]
\[H^{2,0}(X) = H^{2,0}(Y \times E)^{<\pm 1>} =
H^2(Y)^{<\pm 1>} = 0.\]
 Thus $X$ is indeed a smooth Calabi--Yau
threefold. It is non-rigid, because $T(Y\times
E)$ descends to $X$ and each of its Hodge pieces
$H^{p,q}(T(Y\times E))$ is $1$\/-dimensional.

To compute the monodromy around each singular
fiber, we notice that for a generic fiber $X_t =
f^{-1}(t)$ the Kummer structure gives a canonical
decomposition
\[H^2(X_t,\BQ) =
(H^2(A_t,\BQ) \oplus \BQ A_t[2])^{<\pm 1>}.\]
 Our examples were chosen so that the action
of $\pm 1$ on $A_t[2]= E_{\Gamma,t}[2] \times
E_t[2]$ is trivial.  Moreover, in the K\"unneth
decomposition
\[H^2(A_t) = H^2(E) \oplus (H^1(E)\otimes
 H^1(E_{\Gamma,t})) \oplus H^2(E_{\Gamma,t})\]
the $\pm 1$ action is trivial on the first and
last factors, is trivial on $H^1(E)$ and is
unipotent on $H^1(E_\Gamma,t)$ around each
singular fiber of $f$ (namely the cusps of
$\Gamma$). Thus the monodromy of the fibration
$f$ is unipotent as well.

To compute the Kodaira--Spencer map $\Theta(f)$
for our $f$ we embed it into the Kodaira--Spencer
map for $Y\times E \rightarrow
X(\tilde{\Gamma})$. This map is the cup product
with the Kodaira--Spencer class $\Theta$ which
itself is $\Theta_{Y/X(\tilde{\Gamma})} \otimes
\Theta_{X(\tilde{\Gamma}) \times
E/X(\tilde{\Gamma})}$. Since the Kodaira--Spencer
class of a trivial fibration vanishes, it follows
that $\Theta(f) = 0$.

Our examples all have $6$ singular fibers. Hence
\[\frac{1}{2}\deg \Omega^1_{\BP^1}(\log S) =
\frac{1}{2}  \deg\cO_{\BP^1}(-2+6) = 2.\]

On the other hand, since $X$ is a Calabi--Yau
variety we have
\[\omega_{X/\BP^1} = \omega_X\otimes
(f^*\omega_{\BP^1})^{-1} =
(f^*\omega_{\BP^1})^{-1}.\] Hence
\begin{equation}
\label{deg}
 f_*\omega_{X/\BP^1} =
f_*f^*(\omega_{\BP^1})^{-1}=
(f_*f^*\omega_{\BP^1})^{-1}=
\omega_{\BP^1}^{-1},
\end{equation}
whose degree is $2$ as well, concluding the proof
of Theorem 7. \qed
\medskip

{\bf Remarks 8}: {\sl (1) In the other
cases  in Table 2 the monodromy on the points of
order $2$ of $A_t[2]$ is non-trivial, and the
calculation gives that the monodromy of $f$
around the cusps is not unipotent. 
We know by Remark 3,  the groups \#1, \#5 and \#8 are in the
same $PSL(2,\BR)$-conjugacy class. However, this group theoretic
property does not guarantee isomorphisms of the corresponding
Calabi--Yau threefolds, since the fiber structures are not
preserved.  Similarly, the groups \#4, \#6 and \#9 are $PSL(2,\BR)$-conjugate, but
geometric structures are different (as the fibers over the cusps
are different).  The same applies to Examples \#2 and \#7.
Therefore, Examples \#4, \#7, \#8 and \#9 are not covered by our examples. 
Also we do not know how to construct examples corresponding
to Example \#3 in Table 2, for which
$\tilde{\Gamma} = \Gamma_1(7)$. We also do not
know whether our examples are the only ones.

\smallskip

 (2) It is an interesting exercise to compute
the full $L$\/-series of our examples. The
results are as follows: Let $N_+$ (respectively
$N_-$) be the motive of algebraic cycles on $Y$
invariant (respectively anti-invariant) by $\pm
1$ acting on the elliptic fibrations of $Y$. Let
$n_{\pm}$ be the rank of $N_\pm$. Then $n_+ + n_-
= 20$, and if we let $\chi_{\delta^{\prime}}$ denote the
quadratic character cut by
$\BQ(\sqrt{\delta^{\prime}}\,)$ (not necessarily the same quadratic
field pre-determined by the modular group corresponding to
the surface), then 
$N_+= \BZ(1)^{n^{\prime}_+}\oplus\BZ(\chi_{\delta^{\prime}}(1))^{n^{\prime\prime}_+}$
and 
$N_- = \BZ(1)^{n^{\prime}_-}\oplus \BZ(\chi_\delta^{\prime}(1))^{n^{\prime\prime}_-}$. 
Here $n^{\prime}_{+}$ (resp. $n^{\prime\prime}_+$)
denotes the number of cycles defined over
$\BQ$ (resp. $\BQ(\delta^{\prime})$) and similarly for $n^{\prime}_-$ (resp.
$n^{\prime\prime}_-$). We have $n_{\pm}=
(n^{\prime}+n^{\prime\prime})_{\pm}$. 
Then we have

\begin{eqnarray*}
L(H^0,s) & = & L(\BZ,s) = \zeta(s) \\
L(H^1,s) & = & 1 \\
L(H^2,s) &= & L(H^2({\BP}^1\times{\BP}^1),s)L(\BZ(1),s)L(N_+,s)\\
          &=& \zeta(s-1)^{16}\zeta(s-1)^{1+n_+^{\prime}}
L(\BQ\otimes\chi_{\delta^{\prime}},s-1)^{n_+^{\prime\prime}}\\
L(H^3,s) & = & L(T(Y)\otimes H^1(E),s)L(N_-\times H^1(E),s) \\
\quad    & = & L(g_3\otimes g_2,s)L(E, s-1)^{n_-^{\prime}}
\prod_{\delta^{\prime}}L(E\otimes\chi_{{\delta}^{\prime}},s-1)^{n_-^{\prime\prime}}
\end{eqnarray*}

(The higher cohomologies are determined by
Poincar\'e duality.)}   
\medskip

{\bf Lemma 9\ }: In cases \#1, \#2, \#5,
and \#6 in Table 2, we have 
$n_+=14$ and $n_-=6$ (so $n_+-n_-=2+6=8$). 
Furthermore, we have
\begin{eqnarray*}
(n_+^{\prime},n_+^{\prime\prime})=\begin{cases}
(12,2)\quad &\mbox{for \#1} \\
(14,0)\quad &\mbox{for \#2} \\
(13,1)\quad &\mbox{for \#5,\, \#6}\end{cases}
\end{eqnarray*}
and
\begin{eqnarray*}
(n_-^{\prime},n_-^{\prime\prime})=\begin{cases}
(3,3)\quad &\mbox{for \#1} \\
(6,0)\quad &\mbox{for \#2} \\
(5,1)\quad &\mbox{for \#5,\,\#6}\end{cases}
\end{eqnarray*}
In case \#3, $n_+= 11$ and $n_-=9$. (For
the last case, confer the article of Hulek and Verrill [HV] for more
detailed discussion.)
\medskip

For the computations of $n_+$ and $n_-$, confer Proposition 2.4 
of Hulek and Verrill [HV]. The Proof of Lemma 9 will be given at 
the end of Section 4.

\section{Explicit Formulas}
We shall now give explicit formulas for the
weight $3$ forms $g_Y = g_{3,\Gamma}$ for the
examples in Table 2. We will denote the weight
$3$ form in the $i$\/th case by $h_i$. By Remark
3 it suffices to compute $h_i$ for $i=8,7,3,4$,
and then $h_8(\tau) = h_5(\tau/2) = h_1(\tau/4)
$, $h_7(\tau) = h_2(\tau/2)$, and $h_6(\tau) =
h_9(\tau/2) = h_4(\tau/2 -1/2)$. Two kinds of
formulas suggest themselves for the $h_i$\/'s:
as a product of $\eta$\/-functions or as inverse
Mellin transforms of the Dirichlet series
attached to Hecke characters. The second method
is always possible since the $g_Y$\/'s are CM
forms. In [M] Martin determined which modular
forms on $\Gamma_1(N)$ can be expressed as a
product of $\eta$\/-functions. This applies to
cases \#8, \#7, \#3, and \#4 in Table 2. 
Hence the same is also true for the \#6 and \#9
cases. For cases \#3, \#7, and \#8 the
corresponding spaces of cusp forms of weight $3$
are $1$\/-dimensional, hence the conditions in
[M] are satisfied and Martin gives the
corresponding forms as $h_3=
\eta(q)^3\eta(q^7)^3$ and $h_7=
\eta(q^2)^3\eta(q^6)^3$. The modular form in
case \#1 is classically known to be $h_1=
\eta(q)^6$, which implies $h_8= \eta(q^4)^6$.
Lastly $h_2= \eta(q)^3\eta(q^3)^3$, and $h_5=
\eta(q^2)^6$.
For \#4, we have $h_4(q)=\eta(q)^2\eta(q^2)\eta(q^4)\eta(q^8)^2$.
We will prove the following

\medskip

{\bf Proposition 10}: {\sl Let $\chi_i$
be the Hecke character for which $L(h_i,s) =
L(\chi_i,s)$ (so that the inverse Mellin
transform of $L(\chi_i,s)$ is $h_i$). Let
$a_p(h_i)$ be the $p$\/th Fourier coefficient of
$h_i$, and let $K_i = \BQ(\sqrt{\delta_i})$.
Then  we have the following:

(1) The infinite component of
$\chi_i:\BA_{K_i}^\times \rightarrow \BC$ is
given by $\chi_{i,\infty}(z) = z^{-2}$. Moreover
$\chi_i$ is the unique such Hecke character of
conductor $c_i\cO_{K_i}$, where $c_i = 2,2,1,1
\in \cO_{K_i}$ for $i = 8,7,3,4$ respectively.

(2) For each rational prime $p$ which is prime to
the level of the corresponding $\Gamma$, we have
$a_p(h_i) = 0$ if $p$ is inert in $K_i$.
Otherwise, there are $a$, $b$, which are integers
in case \#8 and half integers in the three other
cases, so that $p = a^2+d_i b^2$, where $d_i = 4,
3, 7, 2$ for $i = 8,7,3,4$ respectively. Then $a$
and $b$ are unique up to signs, and $a_p(h_i) =
(a^2-d_ib^2)/2$.}

(3) The newtype of $h_i$ (see (\ref{new})) is the
character defining $K_i$, namely $p\mapsto
\left(\frac{\delta_i}{p}\right)$.
\medskip

{\bf Proof:\ } See e.g. [L2] for the
generalities (in particular regarding the
$\infty$\/-component of $\chi_i$), as well as the
following formula: the conductors of $\chi_i$ and
of $h_i$ are related by
\[\text{cond}(h_i) = \text{Nm}^{K_i}_{\BQ} 
\text{cond}(\chi_i) \text{Disc}(K_i).\]
 Since the level of $h_i$ is respectively $M=16$,
$12$, $7$, and $8$ in cases \#8, \#7, \#3, and
\#4 of Table 2, we get the asserted value of the
$c_i$\/'s, and since all the fields $K_i$ involved have
class number $1$ we have
\begin{equation}
\label{dec}
\BA_{K_i}^\times = (K_i^\times \times U_i \times
\BC^\times)/\mu(K_i)
\end{equation}
where $U_i$ is the maximal compact subgroup
$\hat{\cO}\vphantom{\cO}_{K_i}^\times$ of the
finite id\`eles of $K_i$, and $\mu(K_i)$ is the
group of roots of unity of $K_i$, acting
diagonally (we view $\BC$ as the infinite
completion of $K_i$). The existence and the
uniqueness of $\chi_i$ are then verified in each
case by a straightforward calculation (compare
[L1]).

For the second part, the vanishing of $a_p(h_i)$
for $p$ inert in $K_i$ is a general property of
CM forms. For a split $p$ (prime to
$\text{cond}(h_i)$\/), write $p = a^2 + d_i b^2 =
\text{Nm}^{K_i}_{\BQ} \pi$. Here $\pi$ is a prime
element of $\cO_{K_i}$, so $a$ and $b$ are half
integers. We verify that, up to multiplying $\pi$
by a unit, we can guarantee that $a$ and $b$ are
integers for $i=8$. In all cases, the $a$\/'s and
the $b$\/'s are unique up to signs. Next one
verifies that $\pi\equiv \pm 1
\pmod{c_i\cO_{K_i}}$. Let $\wp$ be the ideal
generated by $\pi$ (notice that changing the sign
of $b$ replaces $\wp$ by its conjugate). Let
$\tr$ denote the trace from $K_i$ to $\BQ$.  By
the general theory, we have that
\[a_p(h_i) = \tr\chi_\wp(\pi) =
\tr\chi_\infty(\pi)^{-1} = \tr\pi^2 =
2(a^2-d_ib^2),\]
 where the second equality holds since the
finite components of $\chi$ other than $\pi$ are
now $1$.

For the third part we notice that the restriction
of $\chi_i$ to $U_i$ in the
decomposition~(\ref{dec}) above gives a Dirichlet
character $\epsilon'_i$ on $\cO_{K_i}$ of
conductor $c_i\cO_{K_i}$. The newtype Dirichlet
character $\epsilon_i$ (on $\BZ$\/) is then the
product $\chi_{K_i}$ by the restriction of
$\epsilon'_i$ to $\BZ$. However, the conductors
$c_i\cO_{K_i}$ of $\chi_i$ are all $1$ or $2$,
and the only character of $\BZ$ of conductor $1$
or $2$ is trivial. Hence the newtype character of
$h_i$ is $K_i$, concluding the proof of
Proposition 10.
\qed
\medskip

{\bf Defining equations for extremal K3
surfaces}: We now discuss how to determine
defining equations for the extremal K3 surfaces
in Theorem 7.  This problem has been getting a
considerable attention lately, for instance,
Shioda [Sh3] and (independently and by a
different method by Y. Iron [I]) have determined
a defining equation for the semi-stable elliptic
K3 surface with singular fibers of type
$I_1,I_1,I_1,I_1,I_1,I_{19}$ whose existence was
established in Miranda and Persson [MP] (this is
given as $\#1$ in their list). As we shall see,
our examples can be determined by a more
classical method.
\medskip

There are several cases where defining equations
can be found in the literature, i.e., Example \#1
in Table 2 is the classical Jacobi quartic
corresponding to $\Gamma(4)$,
\begin{equation}
\label{sig}
y^2=(1-\sigma^2x^2)(1-x^2\sigma^{-2})
\end{equation}
where $\sigma$ is a parameter for $X(4)$. A
Legendre form is given by
\begin{equation}
\label{leg}
Y^2=X(X-1)(X-\lambda)\quad\text{with
$\lambda=\frac{1}{4}(\sigma+\sigma^{-1})^2$}
\end{equation}
(see e.g. Shioda [Sh2]). One checks that the
singular fibers, all of type $I_4$, occur at the
cusps $\sigma = 0,\infty,\pm 1,\pm\sqrt{-1}$.
Moreover the $j$\/-invariant is given by
\[j = 2^4\frac{(1+14 \sigma^2 + \sigma^4)^3}
{\sigma^4(\sigma^4-1)^4}.\]

For the remaining cases, we can find defining
equations using a method due to Tate. Since we
could not find Tate's method in the literature,
we sketch it here. (Actually, we found out after
completing the paper that there are several papers
dealing with this exact issue, e.g., Kubert [K]
and his arguments were reproduced in
Howe--Lepr\'evost--Poonen [HLP]. Also, the paper
of Billing and Mahler [BM] dealt with the same problem.)
\medskip

{\bf A method of Tate to calculate $E_1(N)$\, }:
Let $Y_1(N)$ be the modular curve, and let
$E_1^0(N)\to Y_1(N)$, with $N\geq 4$, be the
universal family of elliptic curves having a
point (or section) $P=P_N$ of order $N$. Tate's
method gives a defining equation for this family
over $\BZ[1/N]$. We start with the general
Weierstrass equation:
\[E:\quad y^2+a_1xy+a_3y=x^3+a_1x^2+a_4x+a_6.\]
Let $P=(x,y)\in E$ be a rational point and assume
that $P, 2P, 3P\neq 0$. Changing   coordinates,
we may put $P$ at $(x,y)=(0,0)$.  So we may
assume $a_6=0$. Since $P$ does not have order
$2$, the tangent line at $(0,0)$ cannot be the
$y$-axis (i.e., $x=0$), which implies that $a_3$
cannot vanish. We can therefore change
coordinates again to obtain $a_4=0$ and the
equation takes the form:
$y^2+a_1xy+a_3y=x^3+a_2x^2$. By making a
dilation, we furthermore get $a_2=a_3$.
Therefore, $E$ has a Weierstrass equation of the
form:
\[(*)\quad y^2+axy+by=x^3+bx^2\quad
\text{with $b\neq 0$}.\] To get a defining
equation for $E_1^0(N)$,  we need to find the
relations on $a$ and $b$ that arise if $P$ has
order $N$. The coordinates of $P, -P, 2P, -2P$
are easily checked to be
\[P=(0,0),\, -P=(-0,-b),\, 2P=(-b, (a-1)b),\,
-2P=(-b,0).\]
 We will also determine the coordinates of $3P$
and of $4P$. At $-2P$ the tangent line is:
\[y=\frac{b}{1-a}(x+b).\]
Substituting this to the equation (*) to get
\[4P=\left(\frac{b}{1-a}+\frac{b^2}{(1-a)^2},\quad
\frac{b^2}{1-a} (1+\frac{b}{(1-a)^2} +
\frac{1}{1-a})\right).\]
 Likewise, the line $x+y+b=0$ intersects $E$ at
$-P$, $-2P$ and $3P$, giving $3P=(1-a,a-1-b)$. We
will give Weierstrass equations for $E_1(N)$ when
$N=4$, $6$, $8$, or $7$:

\medskip $\boxed{E_1(4)\,\,}$ Here we get $a=1$,
giving the equation
\[y^2+xy+ty = x^3 + tx^2\]
(we replaced $b$ by $t$). Here $X_1(4)$ is the
$t$\/-line. By a direct calculation or from
Shioda's result (see also Remark~5), we see that
the singular fibers are over the three cusps, of
types, $I_1^*$, $I_1$, and $I_4$.

\medskip $\boxed{E_1(6)\,\,}$ The equation
$x(4P)=x(-2P)$ readily gives $b= -(a-1)(a-2)$,
giving us the equation for  $E_1(6)$:
\[E_1(6):\quad y^2+axy-(a-1)(a-2)y =
x^3-(a-1)(a-2)x^2.\] Here $a$ is a parameter
(Hauptmodul) on $X_1(6)$. There are four cusps,
and as before one gets from Remark~5 or by a
direct computation that the types are $I_1$,
$I_2$, $I_3$, and $I_6$, matching the widths of
the cusps given in the third column of Table 2.

\medskip $\boxed{E_1(8)\,\,}$ The equation
$y(4P)=y(-4P)$ is equivalent to $ax(4P)+b =
-2y(4P)$. Expanding, cancelling $b$,  and
clearing denominators gives
\[ab(1-a) + (1-a)^2 =
 -2b\left((1-a)(2-a)+b\right).\]
Substituting $b=k(a-1)$ gives
\[(a,b) = \left(\frac{-2k^2+4k-1}{k},
\quad -2k^2+3k-1\right),\]
 Thus $k$ is a parameter on $X_1(8)$ and
$E_1(8)$ (Example \#4) is given by
\begin{equation}
\label{E18}
y^2+\frac{-2k^2+4k-1}{k}xy +(-2k^2+3k-1)y =
x^3+(-2k^2+3k-1)x^2.
\end{equation}
The fibers above the cusps are found as before to
have types $I_1$, $I_1$, $I_2$, $I_4$, $I_8$, and
$I_8$.
\medskip

$\boxed{E_1(7)\,\,}$ In a similar way one gets
the equation for $E_1(7)$ (Example \#3); the result is (see for
instance, Silverman [S] Example 13.4)
\[y^2+(1+t-t^2)xy+(t^2-t^3)y=x^3+(t^2-t^3)x^2.\]
Three singular fibers have type $I_1$, and three
have type $I_7$.
\medskip

$\boxed{E(6,2)\,\,}$ Returning to our cases, we
now handle Example \#2 in Table 2, corresponding to
$\tilde{\Gamma} = \Gamma_0(3)\cap\Gamma(2)$,
whose associated modular curve is
\[X(6,2) = X_1(6) \times_{X_1(2)} X(2)\]
 by pulling $E_1(6)$ back to $X(6,2)$. To do this
we cannot use Tate's method directly, since the
moduli problems associated to $Y(2)$ and to
$Y_1(2)$ are not representable. However $X(2)$ is
the Legendre $\lambda$\/-line, and any elliptic
curve with $\Gamma(2)$\/-level structure can
always be written in Legendre form
\[y^2 = x(x-1)(x-\lambda) =
x(x^2+(-1-\lambda)x+\lambda).\]
 Likewise, given an elliptic curve $E$ with a
point $P$ of order $2$ (over
$\BZ[\frac{1}{2}]$\/), write $E$ in Weierstrass
form $y^2= x(x^2+cx+d)$ where $P=(0,0)$. This
form is unique up to homothety, and hence $c^2/d$
is a parameter on $X_0(2)$. The natural map $X(2)
\rightarrow X_0(2)$ is therefore given by $u =
\frac{(1+\lambda)^2}{\lambda}$. Hence the fibred
product for $E(6,2)$ above is given by equating
\[\frac{(1+\lambda)^2}{\lambda} =
\frac{(4-3a^2)^2}{16(a-1)^3}.\]
 A computation gives that
\[\xi = \frac{32(a-1)^3}{(4-3a^2)(a-2)^2}\left(
\lambda + 1-\frac{(4-3a^2)^2}{32(a-1)^3}\right)\]
is a parameter on $X(6,2)$ such that the map to
$X_1(6)$ is given by
\[a=\frac{2\xi^2-10}{\xi^2-9}.\]
 Under a base change of ramification index $b$
an $I_a$ fiber pulls back to an $I_{ab}$ fiber
(and an $I_a^*$ fiber pulls back to an $I_{ab}^*$
fiber if $a\geq 1$, and $b$ is odd\/). From this or again by
Remark 5 the fiber types are as expected from
Table 2. --- three $I_6$ and three $I_2$ fibers.
\medskip

{\bf Remark}: {\sl Even though we do
not need the following example here, we mention
that a parameter for $X_0(12)$ can be computed
in the same way via the natural map
$X_0(12)\rightarrow X_0(6) = X_1(6)$. By
pull-back this will give a family of elliptic
curves over $X_0(12)$, which will turn out to be
$\Gamma_0(12)$ case (Example \#7) in Table 2. For this we let
$t$ be the parameter for $X_0(4) = X_1(4)$ as
before, and let $a$ be the parameter from before
on $X_0(6)$.  Then $X_0(12)$ is the fibred
product
\[X_1(4)\times_{X_1(2)} X_1(6).\]
To compute the natural ``forgetting'' maps of
$X_1(4)$ and of $X_1(6)$ to $X_1(2)$ we again
bring both $E_1(4)$ and $E_1(6)$ to the form
$y^2=x(x^2+ax+b)$\/:
\[E_1(4): w^2 = v(v^2+(\frac{1}{4}-2b)v+ b^2)\]
(here we completed the square and set
$w=y+(x+b)/2$ and $v=x+b$).
\[E_1(6): w^2 = v(v^2+\frac{4-3a^2}
 {4}v+(a-1)^3).\]
By the above, the fibred product is given by
equating
\[\frac{(\frac{1}{4}-2b)^2}{b^2} =
 \frac{(4-3a^2)^2}{16(a-1)^3}.\]
Thus $a-1$ is a square, say $a=u^2+1$, where $u$
is a parameter on $X_0(12)$ and the pulled-back
family is given by
\[y^2+(u^2+1)xy-u^2(u^2-1)y =
x^3-u^2(u^2-1)x^2.\]}
 One again routinely verifies that the bad fibers
are as expected. (Notice, however, that the
pull-back of the universal family from $X_1(4)$
has $I_a^*$ fibers!)
\medskip

$\boxed{E(8,2)\,\,}$ Next we handle Example \#5 in
Table 2, whose  associated modular curve is
$X(8,2)$. As was explained in Remark 3, we can
take as a parameter the same $\sigma$ as for
$X(4)$ above. However to get the right family, we
will divide the universal elliptic curve by a
section $s$ of order 2. This changes the type of
the singular fibers $E(4)_c$ at a cusp $c$ from
$I_4$ to $I_8$ if $s$ meets $E(4)_c$ at the same
component as the identity section, and to type
$I_2$ otherwise. The singular fibers obtained in
this way are as expected from Table 2.

To get the new family, recall that if an elliptic
curve is given in Weierstrass form
 \[y^2 = x(x^2+ax+b),\]
then the quotient by the two-torsion point
$(0,0)$ is given by the similar equation
\begin{equation}
\label{two}
y^2 = x(x^2-2ax+a^2-4b).
\end{equation}
In particular, for a curve given in Legendre form
$y^2=x(x-1)(x-\lambda)$ the resulting quotient is
$y^2 = x(x^2+2(1+\lambda)x+(1-\lambda)^2)$.
Applying this to the Legendre form (\ref{leg}) of
the Jacobi quartic gives the quotient family in
the form
\[y^2 = x(x^2+ (2+\frac{1}{2}(
\sigma+\sigma^{-1})^2)x +
\frac{1}{16}(\sigma-\sigma^{-1})^4).\]
 As before one sees that the singular fibers have
types $I_8$, $I_8$, $I_4$, $I_4$, $I_4$, and
$I_2$.
\medskip

$\boxed{E(8;4,1,2)\,\,}$ To handle Example \#6 in
Table 2 we proceed in  the same way, dividing the
family $E_1(8)$ by its section of order $2$. The
cusps of $\Gamma_1(8)$ are $N\backslash G/N$,
where $N$ is the upper unipotent subgroup of
G=$SL(2,\BZ/8\BZ)/<\pm\Id>$. Explicitly the cusps
are $\begin{bmatrix} 0\\1 \end{bmatrix} $,
$\begin{bmatrix} 0\\3 \end{bmatrix} $,
$\begin{bmatrix} 1\\2 \end{bmatrix} $,
$\begin{bmatrix} 1\\4 \end{bmatrix} $,
$\begin{bmatrix} 1\\0 \end{bmatrix} $, and
$\begin{bmatrix} 3\\0 \end{bmatrix} $. The
corresponding widths are 1,1,2,4,8, and 8
respectively. We identify the torsion sections of
the universal family $E_1(8)\ra X_1(8)$ with the
subgroup $\begin{pmatrix} *\\0 \end{pmatrix}\in
(\BZ/8\BZ)^2$. Let $s = \begin{pmatrix} 4\\0
\end{pmatrix} $ be the section of order $2$ of
this elliptic fibration. Then $s$  belongs to the
connected component of the $0$\/-section at a
cusp $\begin{bmatrix} a\\b \end{bmatrix} $, if
and only if it is in the subgroup generated by
$\begin{pmatrix} a\\b \end{pmatrix} $. This
happens at the first four cusps above but not at
the last two.

The Kodaira type of the fibers of the quotient
family $E_1(8)/(s)\ra X_1(8)$ (which is our
$E(\Gamma)$  of Example \#6) are accordingly
multiplied by $2$ at the first four cusps, and
divided  by $2$ at the last two cusps. This
results in fiber types $I_8$, $I_4$, $I_4$,
$I_4$, $I_2$, and $I_2$,  in agreement with Table
2.

Starting with our equation (\ref{E18}) for
$E_1(8)$ we set
\[Y = 8k^3\, y
  \left(y+ \frac{-2k^2+4k-1}{2k}x +
  \frac{-2k^2+3k-1}{2} \right)
\qquad\text{and}\qquad X=4k^2(x-k+k^2).\] A
straightforward computation then gives for
$E_1(8)$ the form
\[Y^2 = X(X^2 +(8k^4-16k^3+16k^2-8k+1)X +
 (2k(k-1))^4).\]
By formula (\ref{two}), the quotient family
$E_1(8)/(s)$ is given by
\[y^2 = x(x^2 -2(8k^4-16k^3+16k^2-8k+1)x +
(8k^2-8k+1)(2k-1)^4).\]
\medskip

{\bf Proof of Lemma 9}: We will use the arguments due to
Klaus Hulek and Matthias Sch\"utt for the calculation of 
$n_+$ and $n_-$. The Galois
action on $NS(Y)$ is computed as follows.
Tensoring with $\BQ$, $NS(Y)\otimes\BQ$ has for
basis the (classes of the) general fiber, the
$0$\/-section and those components of the
singular fibers which do not meet the identity
component (the section). The Galois action clearly preserves
the fiber class and the $0$\/-section. The action of
$\pm 1$ on each fiber of type $I_n$ is given
as follows. 
A fiber $I_n$ contributes $n-1$ to the cohomology.
If we enumerate the components $e_1, e_2,\cdots, e_{n-1}$ cyclically,
then $e_j$ will be sent to $e_{-j}$. If $n$ is even,
$n/2$ cycles, $e_{n/2},\, e_j+e_{-j}\, (1\leq j< n/2)$ are fixed
contributing to $n_{+}$; while $e_j-e_{-j}\,(1\leq j< n/2)$ contributing
to $n_{-}$. 
If $n$ is odd, $(n-1)/2$ cycles are fixed contributing to $n_+$,
and equally $(n-1)/2$ cycles to $n_-$. Further, the fields of
definition of the components $e_j$ will determine $n_{\pm}^{\prime}$
and $n_{\pm}^{\prime\prime}$.

In Example \#1, 
the cusps (singularities) are $t=0, \pm 1,\,\infty$ and $\pm\sqrt{-1}$.
Put $i=\sqrt{-1}$.
Then $N_+$ is spanned by the zero-section, the fiber and the
following divisors: $e_{t,2},\, e_{t,1}+e_{t,3}$ where
$t=0, \pm 1, \pm i,\infty$. When $t\in\BQ$ or $t=\infty$, these
divisors are defined over $\BQ$, giving $10$ divisors out of $14$.
Over $t=\pm i$, we see that $e_{i,2}+e_{-i,2}$ and
$(e_{i,1}+e_{i,3})+(e_{-i,1}+e_{-i,3})$
are fixed by complex conjugation, so that these 
are defined over $\BQ$ contributing to $n^{\prime}_+$. 
Hence, as Galois representations, we get
$$N_+=\BZ(1)^{12}\oplus \BZ(\chi_{i}(1))^2$$
so that $n^{\prime}_+=12$, and $n^{\prime\prime}_+=2$.

On the other hand, the space $N_-$ is simply spanned by 
$e_{t,1}-e_{t,3}$ for the six cusps $t$. Over $t=\pm 1$, both are
defined over $\BQ$, contributing to $n^{\prime}_-$. Over $t=0,\infty$,
$e_{t,1}$ and $e_{t,3}$ are conjugate, so the difference is
not fixed under complex conjugation, so it contributes
to $n^{\prime\prime}_-$. Over $t=\pm i$, we have
two divisors $(e_{i,1}-e_{i,3})\pm (e_{-i,1}-e_{-i,3})$. 
One of these is fixed by complex conjugation, while the
other is not. Thus, as Galois representation, 
$$N_-=\BZ(1)^3\oplus \BZ(\chi_{i}(1))^3$$
and reading off the ranks, we get $n^{\prime}_-=3$ and
$n^{\prime\prime}_-=3$.

In Example \#2, 
the cusps are all defined over $\BQ$, the
torsion sections meet all the components of the
fibers (this can be seen either from the moduli
viewpoint or from the equations in the previous 
section). Then $N_+$ is spanned by the zero-section, the fiber
and the divisors $e_{t,1}+e_{t,5},\, e_{t,2}+e_{t,4}$
and $e_{t,3}$ for $I_6$ type singular fibers and
$e_{t,1}$ for $I_2$ type singular fibers. Thus
we compute that 
$n_+=3+3+3+1+1+1+1+1=14$. For the space
$N_-$ is spanned by $e_{t,1}-e_{t,5},\, e_{t,2}-e_{t,4}$ for
$I_6$ type singular fibers. Thus we have
$n_-=2+2+2+0+0+0=6.$
Since all divisors are defined over $\BQ$, all these algebraic
cycles are also defined over $\BQ$, and we have
$$N_+=\BZ(1)^{14}\quad\mbox{and}\quad N_-=\BZ(1)^6.$$
(In particular, this implies that $n_{\pm}^{\prime\prime}=0$
in this case.)

For the other two cases, \#5 and \#6, we use the above
argument to compute $n_{\pm}$. In fact, 
for Example \#5 (resp. \#6),
$$n_+=4+4+1+1+1+1+1+1+1=14 \quad(\mbox{resp.}\,4+2+2+2+1+1+1+1=14),$$
and
$$n_-=3+3=6\quad\mbox{for both cases.}$$
Thus $n_+ = 14$ and $n_-=6$. However, for these examples,
not all algebraic cycles are defined over $\BQ$. In fact,
we use the fact that each of these $K3$ surfaces
is realized as a quadratic base change of a rational modular elliptic
surface (see Top and Yui [TY] for detailed argument).

In the case of Example \#5, this surface is obtained as a pull-back of a
rational elliptic modular surface with $4$ singular fibers
of type $I_4$ over $\\infty$, $I_4$ over $0$, $I_2$ over $1$ and $I_2$
over $-1$. All 
cusps of the pull-back 
over $\infty,\, 0$ and $1$ are defined over $\BQ$. However, the two cusps
of the pull-back above $-1$ are defined only over $\BQ(\sqrt{-1})$.
Put $\sqrt{-1}=i$. Then the divisor $e_{i,1}+e_{-i,1}$ is invariant
under complex conjugation, while the divisor $e_{i,1}-e_{-i,1}$ is not.
Thus, we get
$$N_+=\BZ(1)^{13}\oplus\BZ(\chi_i(1))$$
so that $(n_+^{\prime}, n_+^{\prime\prime})=(13,1).$
On the other hand, all algebraic cycles spanning $N_-$
are defined over $\BQ$ so that 
$N_-=\BZ(1)^6$ and $n_-=n_-^{\prime}=6$.

For Example \#6, the cusps are $t=0,\, \infty,\, \pm 1$ and
$\pm\sqrt{2}$. But the pull-back of the components $e_{0,1}$
and $e_{0,3}$ are conjugate over $\BQ(i)$. This gives
$$N_+=\BZ(1)^{13}\oplus \BZ(\chi_{2}(1))$$
so that $n_+^{\prime}=13$ and $n_+^{\prime\prime}=1$. While 
$$N_-=\BZ(1)^5\oplus \BZ(\chi_i(1))$$
so that $n_-^{\prime}=5$ and $n_-{\prime\prime}=1$.

\qed

{\bf Remark.} {\sl
For \#3, the singular fibers are of type $I_7$ 
and $I_1$ ($3$ copies each).  Hulek and Verrill [HV]
compute that $n_+=11$ and $n_-=9$, and show that all cycles are
defined over $\BQ$. This example does not
admit semi-stable fibrations, but still gives rise to
a non-rigid Calabi--Yau threefold defined over $\BQ$,
and one can still look into the modularity question for the
$L$-series associated to the third cohomology group. This is exactly what
is done in the article of Hulek and Verrill [HV] supplementing
this paper.}
\medskip

{\bf Acknowledgments.}  We thank the
Fields Institute at Toronto, Canada, for its
hospitality. Most of the work described above was
done while the authors were members and
participants in the Automorphic Forms Thematic
Program there during the spring term of 2003. The
second author thanks B. van Geemen, I. Nakamura, A. Sebbar, W.
Stein and K. Ueno for helpful conversations and
correspondences.
After the article was posted on arXiv on April 29, 2003,
we have received feedbacks from several colleagues,
F. Beukers, J. Stientra and J. Top.  We thank them for their interest
and comments.
Finally, the authors are indebted to Matthias Sch\"utt and Klaus
Hulek for their numerous correspondences to clarify the discussions
in Remarks 8 and proof of Lemma 9, in particular, the calculation
of $n_+$ and $n_-$ for our examples.
\bigskip

\hfil{\bf Bibliography}\hfil

\hangafter1 [BTZ] Bartolo, E.A., Tokunaga, H.,
and Zhang, D.-Q., {\it Miranda--Persson's problem
on extremal elliptic K3 surfaces}, Pacific J.\
Math.\ {\bf 202} (2002), 37--72.

\hangafter1 [B] Beauville, A., {\it Les familles
stables  de courbes elliptiques sur $\BP^1$
admettant quatre fibres singuli\`eres}, C. R.
Acad. Sc. Paris, {\bf 294} (1982), 657--660.

\hangafter1 [BL] Besser, A., and Livn\'e, R.,
{\it Universal  Kummer families over Shimura
curves}, in preparation.

\hangafter1 [BM] Billing, G., and Mahler, K., {\it On
exceptional points on cubic curves}, J. London Math. Soc. {\bf 15}
(1940), 32--43.

\hangafter1 [D] Deligne, P., {\it Formes
modulaires et repr\'esentations $\ell$--adiques},
S\'em. Bourbaki, f\'ev 1969, exp. {\bf 355},
1--33, Berlin-Heidelberg-New York, springer
(1977).

\hangafter1 [HLP] Howe,E.W., Leopr\'evost, F., and Poonen, B.,
{\it Large torsion subgroups of split Jacobians of curves
of genus two or three}, Forum Math.\ {\bf 12} (2000),
no. 3, 315--364.

\hangafter1
[HV] Hulek, K. and Verrill, H., {\it On the motive of Kummer varieties
associated to $\Gamma_1(7)$ -- Supplement to the paper : {\it
The modularity of certain non-rigid Calabi--Yau threefolds} by R.
Livn\'e and N. Yui).}

\hangafter1 [I] Iron, Y., {\it An explicit
presentation of  a K3 surface that realizes
$[1,1,1,1,1,19]$}, MSc thesis, Hebrew University
of Jerusalem, 2003.

\hangafter1 [JZ] Jost, J. and Zuo, K., {\it
Arakelov  type inequalities for Hodge bundles
over algebraic varieties}, J. Algebraic Geometry
{\bf 11} (2002), 535--546.

\hangafter1 [KM] Katz, N., and Mazur, B., {\it
Arithmetic  Moduli of Elliptic Curves}, Annals of
Math.\ Studies, {\bf 108}, Princeton University
Press, 1985.

\hangafter1 [KS] Kim, Henry H., and Shahidi, F., {\it Functorial products of
$GL_2\times GL_3$ and the symmetric cube for $GL_2$}. With an appendix by
C. Bushnell and G. Henniart, Ann. of Math. (2) {\bf 155} (2002), no.3,
 837--893.

\hangafter1 [K] Kubert, D., {\it Universal bounds on the torsion of
elliptic curves}, Proc. London Math. Soc. (3) {\bf 33} (1976), no. 2, 193--237.

\hangafter1 [L1] Livn\'e, R., {\it On the conductor of mod $\ell$ Galois
representations coming from modular forms}, J. Number Theory {/bf 31} (1989),
no. 2, 133--141.

\hangafter1 [L2] Livn\'e, R., {\it Motivic orthogonal two-dimensional representations
of $\mbox{Gal}(\bar{\Q}/\Q)$}, Israel J. Math. {\bf 92} (1995), no. 1-3, 149--156.
 
\hangafter1 [M] Martin, Y., {\it Multiplicative
$\eta$\/-quotients}, Trans.\ Amer.\ Math.\ Soc.\
{\bf 348} (1996), 4825--4856.

\hangafter1 [MP] Miranda, R., and Persson, U.,
{\it  Configurations of $I_n$ fibers on elliptic
$K3$ surfaces}, Math.\ Z.\  {\bf 201} (1989),
339--361.

\hangafter1
[S] Schoen, C., {\it On fiber products of rational elliptic surfaces
with section}, Math. Z. {\bf 197} (1988), 177--199.

\hangafter1 [Se] Sebbar, A., {\it Classification
of  torsion-free genus zero congruence groups},
Proc.\ Amer.\ Math.\ Soc.\ {\bf 129}, No. 9
(2001), 2517--2527.

\hangafter1 [SY] Saito, M-H, and Yui, N., {\it
The  modularity conjecture for rigid Calabi--Yau
threefolds over $\Q$}, Kyoto J.\ Math.\ {\bf 41},
no. 2 (2001), 403--419.

\hangafter1
[SZ] Shimada, I., and Zhang, De-Qi, {\it Classification
of extremal elliptic $K3$ surfaces and fundamental groups
of open $K3$ surfaces},
Nagoya Math. J. {\bf 161} (2001), 23--54.

\hangafter1 [Sh1] Shioda, T., {\it Elliptic
modular  surfaces}, J.\ Math.\ Soc.\ Japan {\bf
24}, no. 1 (1972), 20--59.

\hangafter1 [Sh2] Shioda, T., {\it On rational
points  of the generic elliptic curve with level
$N$ structure over the field of modular functions
of level $N$}, J.\ Math.\ Soc.\ Japan {\bf 25}
(1973), 144--157.

\hangafter1 [Sh3] Shioda, T., {\it Discrete
moduli and  integral points}, in preparation.

\hangafter1 [S] Silverman, J., {\it The
Arithmetic of  Elliptic Curves}, Graduate Text in
Mathematics {\bf 106}, Springer--Verlag, New York
1986.

\hangafter1
[STZ] Sun, X., Tan, S.-L. and Zuo, K.,
{\it Families of $K3$ surfaces over curves satisfying the equality
of Arakelov--Yau's type and modularity}, Math. Res. Lett. {\bf 10} (2003), no.2-3,
323--342.

\hangafter1
[TY] Top, J., and Yui, N., {\it Explicit equations of some elliptic
modular surfaces}, Rocky Mountain J. of Math. (to appear).

\hangafter1 [W] Wiles, A., {\it Modular elliptic
curves  and Fermat's last theorem}, Ann.\ of
Math.\ {\bf 141} (1995), 443--551. Taylor, R.,
and Wiles, A., {\it Ring-theoretic properties of
certain Hecke algebras}, Ann.\ of Math.\ {\bf
141} (1995), 553--572.

\hangafter1 [Y] Yui, N., {\it Update on the
modularity of Calabi--Yau varieties}, with appendix by
Verrill, H., {\it The $L$-series of rigid Calabi--Yau threefolds
from fiber products of elliptic curves}, in  {\it Calabi--Yau Varieties and
Mirror Symmetry}, Fields Inst. Commun. {\bf 38} (2001), 307--362, 
Amer. Math. Soc.
\end{document}